

\magnification\magstep1
\baselineskip=14pt
\vsize24.0truecm
\font\csc=cmcsc10
\def\hatt{\widehat}
\def\tilda{\widetilde}
\def\dell{\partial}
\hyphenation{para-meter}
\font\bigbf=cmbx12 

\def\tr{{\rm t}}
\def\dd{{\rm d}}
\def\midd{\,|\,}
\def\E{{\rm E}}
\def\WL{{\rm WL}}

\centerline{\bigbf Estimating the logistic regression equation}

\medskip
\centerline{\bigbf when the model is incorrect}

\bigskip
\centerline{\bf Nils Lid Hjort}

\medskip
\centerline{\sl Norwegian Computing Centre and University of Oslo}

\medskip
\centerline{\sl -- January 1990 --}

\bigskip
{{\narrower\narrower\smallskip\baselineskip12pt\noindent
{\csc Abstract.} Protesting mildly against the notion of an exactly correct
parametric model the view is adopted that the logistic regression equation
is merely an approximation to the underlying, true function. The behaviour 
of likelihood based estimators is investigated in such a general framework.
The maximum likelihood estimator is shown to be consistent for a certain
least false parameter value minimising a weighted average of
quantities that measure the distance from the true to the parametric model.
Asymptotic normality is also demonstrated. Finally a number of additional
remarks are offered, 
some pointing to natural generalisations and some to new questions
for research, like weighted and local likelihood estimation methods.

\smallskip\noindent 
{\csc Key Words:} {\sl incorrect model; 
agnostic parameter estimation; local likelihood;
logistic regression} \smallskip}}

\bigskip
\centerline{\bf 1. Introduction}

\medskip\noindent  
Let $q$ be the probability of some particular event of interest for an 
individual object or `case'. Assume that variables $x_1,\ldots,x_d$ can
be measured for each individual case, and that these 
are thought to influence the
probability $q$, so that $q = q(x_1,\ldots,x_d)$. A common way of modelling
this is via the logistic regression equation
$$q= {{\exp(\alpha+\beta^\tr x)}\over{1 + \exp(\alpha+\beta^\tr x)}}.\eqno(1.1)$$
The interpretation is that if the relative frequency of the event in question 
were observed for a large number of cases with identical covariates
$x = (x_1,\ldots,x_d)$, then this relative frequency would be close to the
right hand side of (1.1), for a particular set of parameters 
$(\alpha,\beta_1,\ldots,\beta_d)$. 

The parameters are usually estimated using the likelihood resulting from 
observing whether or not the event occurs for a number of cases with
recorded covariates, i.e.
$$L = \prod_{i=1}^n q_i^{z_i}(1-q_i)^{1-z_i},\eqno(1.2)$$
where $z_i$ is one or zero depending upon whether the event occurs or not for
case no. $i$, and $q_i = q_\beta(x_i) = q_\beta(x_{i,1},\ldots,x_{i,d})$. 
The theory that supports traditional inference is asymptotic and with the firm 
underlying assumption that the model specification (1.1) in fact is correct.

In the present note the behaviour of likelihood-based estimation is studied
under the more realistic and less restrictive assumption that $q$ simply is
some fixed, unknown function of $(x_1,\ldots,x_d)$. It it not so clear what
the maximum likelihood estimators take aim at now when the model is 
incorrect; this is clarified in Section 2. Section 3 provides the limit 
distribution of the maximum likelihood 
estimator, generalising the classical model-based
result. Finally some relevant additional remarks are offered in Section 4.

\bigskip
\centerline{\bf 2. What is the maximum likelihood 
				estimator really estimating?}

\medskip\noindent 
It is notationally convenient to write
$\alpha + \sum_{j=1}^d\beta_jx_j =\beta^\tr x$ now,
where $\beta = (\beta_0,\ldots,\beta_d)^\tr$ with $\beta_0 = \alpha$,
and $x = (1,x_1,\ldots,x_d)^\tr$. Consider
$$\eqalign{{1\over n} \log L &= {1\over n}\sum_{i=1}^n 
	\bigl\{z_i\log q_i + (1-z_i)\log(1-q_i)\bigr\} \cr
	&= {1\over n}\sum_{i=1}^n
	\bigl[z_i\beta^\tr x_i - \log\{1+\exp(\beta^\tr x_i)\}\bigr], \cr}\eqno(2.1)$$
and let $\hatt\beta = (\hatt\beta_0,\ldots,\hatt\beta_d)^\tr$ be the 
maximum likelihood estimators,
i.e.~the parameter values that maximise the above expression. We shall study
the limit distribution properties of $\hatt\beta$ without assuming (1.1) to
be exact.

Suppose that the empirical distribution of $x_1,\ldots,x_n$ becomes close to
some $H(\dd x)$, the `covariate distribution', in a sense to be made precise
later on. Sometimes the covariates are
under the experimenter's control, and in such cases the sampling design
employed of course influences $H$. But in many cases the $x_i$'s are rather
drawn randomly from some population of covariate vectors, and then $H$ is
simply the probability distribution of such a randomly chosen $x_i$. See
also the remarks made in Section 4.

Under regularity conditions ${1\over n} \log L(\beta)$ converges in probability
to
$$\eqalign{\E\,\bigl[z\beta^\tr x - \log\{1+\exp(\beta^\tr x)\}\bigr]
	&= \E\,\bigl[q(x)\beta^\tr x - \log\{1+\exp(\beta^\tr x)\}\bigr] \cr
	&= \int\bigl[q(x)\beta^\tr x - 
		\log\{1+\exp(\beta^\tr x)\}\bigr]\,H(\dd x),\cr}\eqno(2.2)$$
for each fixed $\beta$, in which $q(x)$ is the underlying {\sl true} 
probability, i.e.~$z_i\midd x_i$ is Bernoulli $\{1,q(x_i)\}$. 
Here and in the following `$\E$' denotes expectation w.r.t.~the 
simultaneous distribution of $(x,z)$, i.e. $x \sim H(\dd x)$ and 
$\Pr\{z = 1\midd x\} = q(x)$. --- The maximum likelihood estimator 
therefore in effect takes aim at the parameter value $\beta^0 =
(\beta_0^0,\ldots,\beta_d^0)^\tr$ that maximises this expression. This is a
sensible statistical aim even if the parametric model is incorrect, in that
$\beta^0$ minimises the following natural distance measure from the true
$q(.)$ to the parametric $q_\beta(.)$:
$$\Delta\{q(.),q_\beta(.)\} = 
		\int D\{q(x),q_\beta(x)\}\,H(\dd x),\eqno(2.3)$$
$$D\{q(x),q_\beta(x)\} = q(x)\log {q(x)\over q_\beta(x)}
	+\{1-q(x)\}\log {{1-q(x)}\over{1-q_\beta(x)}}.\eqno(2.4)$$
Thus $D\{q(x),q_\beta(x)\}$ is the Kullback-Leibler information-theoretic
distance from Bernoulli $\{1,q(x)\}$ to Bernoulli $\{1,q_\beta(x)\}$, and
$\Delta\{q,q_\beta\}$ is a weighted average of these $x$-conditional
distances. One has $\Delta = 0$ only if there is a $\beta^1$ with 
$q(x) = q_{\beta^1}(x)$ for all $x$ in the support of $H$, i.e.~only if
the logistic regression equation (1.1) is exactly fulfilled (for the
$x$-region of interest), and then $\beta^0 = \beta^1$.

Let us in general call $\beta^0$ minimising $\Delta$ above the {\sl least
false} parameter value, a small but statistically important conceptual leap
from the classical notion of a `true' parameter value.
It is shown in the following section that $\hatt\beta$
is consistent for this least false parameter value. 

\bigskip
\centerline{\bf 3. Limit distribution of the maximum likelihood estimator}

\medskip\noindent 
We get from (2.1) and (1.1) that
$${1\over n}\dell\log L/\dell\beta_u = {1\over n}\sum_{i=1}^n 
	x_{iu}\{z_{i}-q_\beta(x_i)\}$$
and that
$${1\over n}\dell^2\log L/\dell\beta_u\beta_v = 
	-{1\over n}\sum_{i=1}^n x_{iu}x_{iv}\,q_\beta(x_i)\{1-q_\beta(x_i)\},$$
so that the Fisherian information matrix becomes $n$ times
$$J_n = J_n(\beta) = -{1\over n}\dell^2\log L/\dell\beta\dell\beta
	= {1\over n}\sum_{i=1}^n x_ix_i^\tr\,q_i(1-q_i),\eqno(3.1)$$
again writing $q_i = q_\beta(x_i)$ for short. 
The traditionally employed asymptotic
result on $\hatt\beta$ is that $\hatt\beta \approx N_{d+1}\{\beta,
\hbox{$1\over n$}J_n(\hatt\beta)^{-1}\}$, viz. statistical computer program 
packages like {\csc bmdp, sas} etc. However, $\sqrt{n}(\hatt\beta-\beta^0)$
does {\sl not} converge in distribution to $N_{d+1}\{0,J(\beta^0)^{-1}\}$, 
where $J(\beta^0)$ is the limit in probability of $J_n(\hatt\beta)$, 
under the present agnostic state of affairs.

\smallskip
{\csc Theorem.} Suppose that the covariates $x_1,\ldots,x_n$ can be assumed to
be independently drawn from a `covariate distribution' $H(\dd x)$. If 
$\int\Vert x\Vert\,H(\dd x) < \infty$, then there is a unique least false 
parameter $\beta^0$ minimising the distance measure (2.3), and $\hatt\beta$
converges in probability to $\beta^0$. If $\int xx'\,H(\dd x)$ is 
finite too, then 
$$\sqrt{n}(\hatt\beta-\beta^0)\rightarrow_d 
		N\{0,J(\beta^0)^{-1}K(\beta^0)J(\beta^0)^{-1}\},$$
in which 
$$J(\beta) = \E\,xx^\tr \,q_\beta(x)\{1-q_\beta(x)\} 
	= \int xx^\tr q_\beta(x)\{1-q_\beta(x)\}\,H(\dd x),$$
$$K(\beta) = \E\,xx^\tr \{z-q_\beta(x)\}^2 
	= \int xx^\tr\bigl[\{1-q_\beta(x)\}^2q(x) 
			+ q_\beta(x)^2\{1-q(x)\}\bigr]\,H(\dd x).$$
Consistent estimators for $J(\beta^0)$ and $K(\beta^0)$ are respectively
$$\hatt J = J_n(\hatt\beta) 
	= {1\over n}\sum_{i=1}^n x_ix_i^\tr\hatt q_i(1-\hatt q_i), \quad 
	\hatt K = {1\over n}\sum_{i=1}^n x_ix_i^\tr(z_i-\hatt q_i)^2,$$
in which $\hatt q_i = q_{\hatt\beta}(x_i) = 
\exp(\hatt\beta^\tr x_i)/\{1+\exp(\hatt\beta^\tr x_i)\}$. 

\smallskip
{\csc Proof:} The Fisher information matrix $J_n(\beta)$ is 
everywhere positive definite,
which implies that the log-likelihood function is concave in $\beta$. This
is pleasant seen from the point of view of numerical optimisation, i.e.
$\hatt\beta$ is unique and easy to find, 
and even Bayes estimators, expressed as
$d+1$-dimensional integrals, can be computed stably. But the log-concavity
is also pleasant from a theoretical point of view: the pointwise limit (2.2)
in probability of $(1/n) \log L$ is also concave, and the maximiser of 
$(1/n)\log L$ is consistent for the maximiser of the limit function,
cf. Andersen and Gill (1982, Appendix 2). This 
proves $\hatt\beta\rightarrow_p\beta^0$. (The moment condition 
$\int |x_u|\,H(\dd x) < \infty$ for $u = 1,\ldots,d$ is what is needed to ensure
existence of and pointwise convergence in probability to the (2.2) function.)

The least false $\beta^0$ minimises the $\Delta$-distance w.r.t. $\beta$, 
and has accordingly the property that the partial derivatives vanish at that
point, i.e.
$$\int x_u\{q(x)-q_{\beta^0}(x)\}\,H(\dd x) = 0,\quad u = 0,\ldots,d,$$
which is another way of quantifying the way in which $\beta^0$ makes
$q_\beta(x)$ come as close as possible to the true $q(x)$. 
Hence $(1/n)\dell\log L(\beta^0)/\dell\beta\rightarrow_p 0$, and
$${1\over\sqrt{n}} {{\dell\log L(\beta^0)}\over{\dell\beta}}
	= {1\over\sqrt{n}}\sum_{i=1}^nx_i\{z_i-q_{\beta^0}(x_i)\}
	\rightarrow_d N_{d+1}\{0,K(\beta^0)\}$$
by the central limit theorem. The standard Taylor argument gives
$$\sqrt{n}(\hatt\beta-\beta^0) = J_n(\tilda\beta)^{-1} 
	{1\over\sqrt{n}} {{\dell\log L(\beta^0)}\over{\dell\beta}},$$
in which $\tilda\beta$ is somewhere between $\beta^0$ and $\hatt\beta$,
and the limit distribution result is proved as soon as it is demonstrated
that $J_n(\tilda\beta)\rightarrow_p J(\beta^0)$. But $J_n(\beta^0)$ goes
a.s. to $J(\beta^0)$ by the law of large numbers, and 
$J_n(\beta_n)-J_n(\beta^0)$ can be shown to go to zero in probability 
whenever $\beta_n\rightarrow\beta^0$, proving the assertion, because
$\hatt\beta$ is consistent. The latter arguments also prove that $\hatt J$
is consistent for $J(\beta^0)$, and can be parallelled to give
$\hatt K \rightarrow_p K(\beta^0)$. $\heartsuit$

\bigskip
\centerline{\bf 4. Concluding remarks}

\medskip\noindent
Some concluding comments are as follows, some pointing
to further questions and research. 

\smallskip 
{\csc Remark A.} The theorem shows that inference about the logistic 
regression coefficients $\beta_u$ can be carried out using 
$\hatt\beta
\approx N_{d+1}\{\beta^0,{1\over n}\hatt J^{-1}\hatt K \hatt J^{-1}\}$.
This constitutes a basically more sound and robust alternative to the 
traditional one, in which $\hatt\beta \approx 
N_{d+1}\{\beta,{1\over n}\hatt J^{-1}\}$ is used.

\smallskip
{\csc Remark B.} If $q(x) \equiv q_{\beta^0}(x)$ then 
$J(\beta^0) = K(\beta^0)$, and only then is the classic 
$\sqrt{n}(\hatt\beta-\beta^0)\rightarrow_d N_{d+1}\{0,J(\beta^0)^{-1}\}$
in force. One may indeed base a goodness of fit test for the parametric
model on how close $\hatt J$ is to $\hatt K$.

\smallskip
{\csc Remark C.} If $p(\beta)\,d\beta$ is some prior distribution density
for $\beta$ with finite expectation, then the Bayes solution 
$$\beta^* = {{\int\beta L(\beta)p(\beta)\,\dd\beta}\over
		{\int L(\beta)p(\beta)\,\dd\beta}}$$
can be shown to be asymptotically equivalent to the maximum likelihood
estimator, i.e. $\sqrt{n}(\hatt\beta-\beta^*)\rightarrow_p 0$, regardless
of the actual prior density, regardless of whether the logistic
regression model $q(x) \equiv q_\beta(x)$ holds or not, and regardless
of the covariate distribution $H(\dd x)$. In particular 
$\sqrt{n}(\beta^*-\beta^0)$
has the same limit distribution as the one given in the theorem for the
maximum likelihood eastimator, and provides a frequentist justification for
the Bayes solution. These claims can be proved using arguments similar to
those of Hjort (1986a, Section 2).

\smallskip
{\csc Remark D.} The `covariate distribution' $H(\dd x)$ is the conceptual
limit of the empirical distribution of the actually employed $x_1,\ldots,x_n$.
If these are under the experimenter's control and are chosen according to
some systematic sampling design (say at points uniformly spread 
in a rectangle in the $x$-space), 
then $H$ is the corresponding limit distribution for this
sampling strategy (say the continuous uniform measure on the appropriate
$x$-rectangle). In particular it should be noted that the exact value of 
the least false parameter $\beta^0$ depends upon the $H$-distribution 
(except in the unrealistic ideal case where the logistic regression 
equation holds exactly; this ideal situation presumably prevails 
only in statisticians'
computer simulations). Accordingly, if two laboratories study the same
phenomenon, with the very same $q(x) = \Pr\{z = 1\midd x\}$, but adopt different
sampling procedures for the chosen $x$'s, then their two published 
$\beta$-estimates will in fact aim at two different parameter values.

\smallskip
{\csc Remark E.} Several other published methods of estimating $\beta$ 
may all be consistent for the `true' $\beta$ in the classic framework
of a true model (1.1), but they may easily be consistent for different
least false parameters when (1.1) does not hold. This serves to point
out that different approaches to estimation correspond (in general) to
different modes of `least false'-ness. (Remark C contains the pleasing
exception to the rule; maximum likelihood and Bayes estimators have
the same aims and the same asymptotic behaviour even when the model on
which the likelihood is based is incorrect.)

The notion of an (asymptotically) optimal estimator becomes less clear
in view of this. In particular, a theorem proving the optimality of
the maximum likelihood under exact model conditions appears perhaps
less impressive and less interesting.

\smallskip
{\csc Remark F.} The Fisher information notion is also less clear in the
absence of the luxurious support of a `true model'. Perhaps 
$(J^{-1}KJ^{-1})^{-1} = JK^{-1}J$ ought to take $J$'s place.

\smallskip
{\csc Remark G.} It is perhaps disturbing that the distribution of 
$x$-variables influences the most fitting parameter value. This can however
also be put to the statistician's advantage, for example by estimating 
$\beta$ separately in different $x$-regions of interest. A general
`weighted likelihood' generalisation of (1.2) is 
$$\WL = \prod_{i=1}^n\bigl[q_\beta(x_i)^{z_i}
			\{1-q_\beta(x_i)\}^{1-z_i}\bigr]^{w(x_i)}$$
for an appropriate weight function $w(x)$, for example one that is one in
a particular $x$-region and zero elsewhere. The corresponding
`maximum weighted likelihood estimator' maximises
$${1\over n}\log \WL = {1\over n}\sum_{i=1}^n
 \bigl[z_i\beta^\tr x_i-\log\{1+\exp(\beta^\tr x_i)\}\bigr]\,w(x_i),\eqno(4.1)$$
and converges to a new least false parameter value that minimises
$$\Delta_w\{q,q_\beta\} = \int D\{q(x),q_\beta(x)\}w(x)\,H(\dd x),$$
cf. (2.3) and (2.4).

\smallskip
{\csc Remark H.} In many real world situations the $x$'s can be 
assumed to constitute measurement vectors for some randomly sampled 
underlying population. 
In such cases
$H$ has a mixture density 
$$h(x) = \pi_0 f_0(x) + \pi_1 f_1(x),\eqno(4.2)$$
where $f_0$ and $f_1$ are the two class-distributions involved;
the distribution of $x$'s given that $z=0$ and given that $z=1$, respectively.
It also follows that
$$q(x) = {{\pi_1 f_1(x)}\over{h(x)}} = 
	{{\pi_1 f_1(x)}\over{\pi_0 f_0(x) + \pi_1 f_1(x)}},\eqno(4.3)$$
which can be substituted in equations above for $\Delta\{q,q_\beta\}$, 
$J(\beta)$, and $K(\beta)$, to give additional understanding of the quantities
involved.

\smallskip
{\csc Remark I.} If the sampling situation underlying (4.1) and (4.2) 
prevails, then the logistic regression equation amounts to 
$$\log{{\pi_1 f_1(x)}\over{\pi_0 f_0(x)}} = 
			\beta_0 + \sum_{j=1}^d\beta_jx_j,$$
i.e.~a modelling of the ratio of the two class densities. Sometimes more
explicit parametric models are employed that actually entail such a 
log-linear relationship. For example, if $(x_1,x_2)$ is distributed as
${\rm Beta}(a_0,b_0)\,\times\,{\rm Beta}(c_0,d_0)$ for
$z=0$ cases and as 
${\rm Beta}(a_1,b_1)\,\times\,{\rm Beta}(c_1,d_1)$ for
$z=1$ cases, then the log-ratio is linear in $\log x_1$, $\log(1-x_1)$,
$\log x_2$, $\log(1-x_2)$, with corresponding logistic regression 
parameters $\beta_1 = a_1-a_0$, $\beta_2 = b_1-b_0$, $\beta_3 = c_1-c_0$, 
$\beta_4 = d_1-d_0$. Similarly and more familiarly, if 
$x\sim N_d\{\mu_0,\Sigma_0\}$ and $x\sim N_d\{\mu_1,\Sigma_1\}$
for the two classes, then $q(x)$ is linear 
in $x_i$'s and $x_ix_j$'s. 

In such situations 
an alternative logistic regression estimation strategy is to 
carry out parameter estimation separately for the 
two groups, leading for example to $\hatt\beta_1 = \hatt a_1-\hatt a_0$ etc.
in the Beta-example. If the model is right then this second strategy is
more efficient, because it utilises more fully the information
available, in a Fisherian sense. A fortiori the resulting discriminant
analysis rule, assigning `class 1' to incoming $x$'s whenever
$q(x;\hatt\beta) > \hbox{$1\over2$}$ and `class 0' otherwise, is better
with $\hatt\beta$ arrived at by group-wise estimation.
Explicit calculations are carried
out by Efron (1975) in the Gaussian case and in unpublished work 
by the present author in the Beta situation. The logistic regression approach
has however an edge regarding robustness, since a whole class of models
for the class densities $f_0$ and $f_1$ lead to the same log-ratio.

\smallskip
{\csc Remark J.} Suppose that parametric models are used for the class
densities and that these are consistent with (1.1) for the log-ratio.
Subscribing to the same view of parametric modelling as elsewhere 
in this paper 
we are led to consider them as mere approximations to the complex reality.
As a consequence the parameter estimates stemming from group-wise 
parametric estimation also aim at certain least false parameters, 
giving in their turn new least false parameters for the logistic
regression problem. These are not the same as those discussed in
earlier sections, and correspond to a different distance measure than (2.3).
See also Remark E, and Hjort (1986b, Chapter 3).

\smallskip
{\csc Remark K.} The reasoning in this paper applies equally well 
to models specifying other
parametric structures than (1.1) for $q(x)$, for example in probit analysis,
where $q(x) = \Phi(\beta^\tr x)$ and $\Phi$ is the standard Gaussian cumulative
distribution function. One can again prove that the maximum likelihood
estimator is consistent for a particular least false parameter value, 
and that it has a Gaussian limit distribution, with a 
limit distribution covariance matrix of the $J^{-1}KJ^{-1}$ structure. 

\smallskip
{\csc Remark L.} The line of arguments is in fact even more general. The
reasonably general result about maximum likelihood estimation in parametric
models is that the estimator is consistent for the parameter value that
is least false according to the Kullback-Leibler distance between the
true model and the parametric approximation to the model, and that it
is asymptotically normal with covariance matrix of the $J^{-1}KJ^{-1}$
type, see Hjort (1986b, Chapter 3) for more details.

\smallskip
{\csc Remark M.} Let us finally mention nonparametric
alternatives to the traditional logistic regression. Work by 
Hastie and Tibshirani (1986, 1987), Green (1987), and Azzalini,
Bowman, and H\"ardle (1988) discuss different smoothing techniques
for the direct estimation of $q(x) = \Pr\{z=1\midd x\}$. 

Here is another one,
which can be computationally involved but appears natural, and which can
be considered as an extreme version of the weighted likelihood method
of Remark G. Fix a covariate value $x_0$, and let $w(\cdot)$ be one on a
neighbourhood $N(x_0)$ of $x_0$ and zero elsewhere. Obtain the 
maximum weighted 
likelihood estimator, say $\hatt\beta(x_0)$, and in the end use
$$\tilda q(x) = {{\exp\{\hatt\beta(x_0)^\tr x\}}\over
			{1+\exp\{\hatt\beta(x_0)^\tr x\}}}$$
for $x$ close to $x_0$ (and perhaps only for $x = x_0$).
This corresponds to using $\hatt a(x_0) + \hatt b(x_0)(x-x_0)$, $x$ 
close to $x_0$, as a nonparametric estimate of $\E\,(Y\midd x)$ in 
ordinary regression, with locally obtained $\hatt a(x_0)$ and 
$\hatt b(x_0)$. --- Many similar suggestions can be made. One of them is
to use 
$$q^*(x) = {{\exp\{\hatt\beta(x)^\tr x\}}\over
			{1 + \exp\{\hatt\beta(x)^\tr x\}}},$$
where $\hatt\beta(x)$ maximises $\sum_{i=1}^n [z_i\beta^\tr x_i - 
\log\{1+\exp(\beta^\tr x)\}]\,K(x-x_i)$ for some kernel 
weight function $K(\cdot)$.

\smallskip
{\csc Remark N.} A simpler and presumably
better solution in the case of randomly samp\-led $x$'s is to use
$$\hatt q(x) = {{\pi_1 \hatt f_1(x)}\over
	{\pi_0 \hatt f_0(x) + \pi_1 \hatt f_1(x)}},$$
in which $\hatt f_0$ and $\hatt f_1$ are nonparametric density estimates.
Such a procedure, involving separate and interesting problems regarding
choice of (perhaps two different) smoothing parameters, 
has apparently not been studied in the literature. 



\bigskip
\centerline{\bf References}

\medskip 
\item{1} Andersen, P.K.~and Gill, R. (1982). 
Cox's regression model for counting processes: a large sample study. 
{\sl Ann.~Statist.~\bf 10}, 1100--1120.

\item{2} Azzalini, A., Bowman, A.W., and H\"ardle, W. (1988). 
On the use of nonparametric regression for model checking. 
{\sl Biometrika~\bf 75}, to appear.

\item{3} Efron, B. (1975). 
The efficiency of logistic regression compared
to normal discriminant analysis. 
{\sl J.~Amer.~Statist.~Assoc.~\bf 70}, 892--898.

\item{4} Green, P.J. (1987). 
Penalized likelihood for general semi-parametric
regression models. 
{\sl Int.~Statist.~Rev.~\bf 55}, 245--259.

\item{5} Hastie, T.J.~and Tibshirani, R. (1986). 
Generalized additive models (with discussion). 
{\sl Statist.~Science~\bf 1}, 297--318.

\item{6} Hastie, T.J.~and Tibshirani, R. (1987). 
Non-parametric, logistic and proportional odds regression. 
{\sl Appl.~Statist.~\bf 36}, 260--276.

\item{7} Hjort, N.L. (1986a). 
Bayes estimators and asymptotic efficiency in
parametric counting process models. 
{\sl Scand.~J.~Statist.~\bf 13}, 63--85.

\item{8} Hjort, N.L. (1986b). 
{\sl Theory of Statistical Symbol Recognition.} 
Norwegian Computing Centre Report No. 778.

%
%
%
%

\bye